\documentclass[letterpaper, 10 pt, conference, twoside]{IEEEtran}
\usepackage{fancyhdr}
\usepackage{url}
\usepackage{amsmath} 
\usepackage{amssymb}  
\usepackage{amsthm}
\usepackage{lipsum,adjustbox}
\usepackage{tikz}
\usepackage{tikz-qtree}
\usepackage[font=footnotesize]{caption}
\usepackage{multirow}

\usepackage[shortlabels]{enumitem}
\renewcommand{\theenumi}{(\roman{enumi})}

\usepackage{dsfont}
\usepackage{cite}
\usepackage{subfigure}
\usetikzlibrary{arrows,shapes,positioning,shadows,trees}

\DeclareMathOperator{\nulls}{\operatorname{null}}
\DeclareMathOperator{\lie}{\mathcal{L}}
\DeclareMathOperator{\Lap}{\mathsf{L}}
\DeclareMathOperator{\A}{\mathsf{A}}
\DeclareMathOperator{\D}{\mathsf{D}}
\DeclareMathOperator{\V}{\mathcal{V}}
\DeclareMathOperator{\E}{\mathcal{E}}
\DeclareMathOperator{\G}{\mathcal{G}}
\newcommand{\QQ}{\mathbf{Q}}
\newcommand{\ww}{\mathbf{w}}
\newcommand{\yy}{\mathbf{y}}
\newcommand{\ee}{\mathbf{e}}

\newcommand{\PiPi}{\mathbf{\Pi}}
\newcommand{\real}{\mathbb{R}}
\newcommand{\realpos}{\mathbb{R}_{>0}}

\newcommand{\integer}{\mathbb{Z}}

\newcommand{\diag}{\operatorname{diag}}
\DeclareMathOperator{\one}{\mathbf{1}}
\DeclareMathOperator{\onem}{\mathds{1}}
\DeclareMathOperator{\zero}{\mathbf{0}}
\DeclareMathOperator{\xx}{\mathbf{x}}
\DeclareMathOperator{\bb}{\mathbf{b}}
\DeclareMathOperator{\LL}{\mathbf{L}}
\DeclareMathOperator{\AAA}{\mathbf{A}}

\DeclareMathOperator{\PP}{\mathbf{P}}
\DeclareMathOperator{\zz}{\mathbf{z}}

\DeclareMathOperator{\grad}{\psi_{\textnormal{grad}}}
\DeclareMathOperator{\gdac}{\psi_{\textnormal{gdac}}}

\DeclareMathOperator{\s}{\textnormal{ss}}
\DeclareMathOperator{\dist}{\psi_{\textnormal{dist}}}
\newcommand{\vv}{\mathbf{v}}
\newcommand{\VV}{\mathbf{V}}
\newcommand{\dd}{\mathbf{d}}

\newtheorem{theorem}{Theorem}[section]
\newtheorem{proposition}[theorem]{Proposition}
\newtheorem{lemma}[theorem]{Lemma}

\theoremstyle{remark}
\newtheorem{remark}{Remark}
\theoremstyle{definition}

\newcommand{\longthmtitle}[1]{\mbox{} \emph{(#1):}}
\newcommand{\until}[1]{\{1,\dots,#1\}}

\newcommand{\oprocendsymbol}{\hbox{$\bullet$}}
\newcommand{\oprocend}{\relax\ifmmode\else\unskip\hfill\fi\oprocendsymbol}

\allowdisplaybreaks

\title{Solving Linear Equations with Separable Problem Data over
  Directed Networks\thanks{This work was partially supported by NSF
    Awards 1917177 and 1947050.}}

\author{Priyank Srivastava \quad Jorge
  Cort\'{e}s  
\thanks{The authors are with the Department of Mechanical and
    Aerospace Engineering, UC San Diego,
    {\tt\small \{psrivast,cortes\}@ucsd.edu}}%
}

\graphicspath{ {Figures/} }

\begin{document}

\maketitle
\thispagestyle{fancy}

\begin{abstract}
  This paper deals with linear algebraic equations where the global
  coefficient matrix and constant vector are given respectively, by
  the summation of the coefficient matrices and constant vectors of
  the individual agents.  Our approach is based on reformulating the
  original problem as an unconstrained optimization.  Based on this
  exact reformulation, we first provide a gradient-based, centralized
  algorithm which serves as a reference for the ensuing design of
  distributed algorithms.  We propose two sets of exponentially stable
  continuous-time distributed algorithms that do not require the
  individual agent matrices to be invertible, and are based on
  estimating non-distributed terms in the centralized algorithm using
  dynamic average consensus.  The first algorithm works for
  time-varying weight-balanced directed networks, and the second
  algorithm works for general directed networks for which the
  communication graphs might not be balanced.  Numerical simulations
  illustrate our results.
\end{abstract}
%\vspace*{-3ex}

\section{Introduction}\label{sec:intro}
\IEEEPARstart{T}{he} importance of solving linear algebraic equations is paramount.
They appear frequently in core mathematics as well as in applications,
in physics and engineering. Nonlinear systems can often be well
understood by their linear approximation.  Due to the recent
development of large-scale networks coupled with parallel processing
power and fast communication capabilities, there is a growing effort
aimed at developing distributed algorithms to solve systems of linear
equations.  Distributed algorithms preserve the privacy of the agents,
are robust against single point of failures, and scale well with the
network size.  Keeping these considerations in mind, this paper is a
contribution to the growing body of distributed algorithms to solve
linear algebraic equations.

\emph{Literature Review:} Justifying the ubiquity of linear equations,
there is a vast and expanding literature to solve them efficiently,
cf.~\cite{DPB-JNT:97,SM-JL-ASM:15,BDOA-SM-ASM-UH:16} and references
therein. However, most of the works consider the information structure
where each agent knows some rows of the coefficient matrix and the
constant vector.  In those cases, the collective problem has a
solution if and only if the individual equations are
solvable. Instead, the problem structure considered here is different,
and assumes that each agent has a full coefficient matrix and constant
vector of its own.  This setting appears frequently in distributed
sensor fusion, where sensors are spatially distributed and they seek
to build a global state estimate (e.g., about the location of a source
or the position of a target) from local measurements,
cf.~\cite{DPS-ROS-RMM:05-ifac,LX-SB-SL:05}.  All the works in this
category rely on the communication graph being \emph{undirected}.  The
work~\cite{DPS-ROS-RMM:05-ifac} relies on the positive definiteness of
the individual matrices to compute the updates and prove
stability.~\cite{LX-SB-SL:05} uses element-wise average consensus for
the coefficient matrix as well as the constant vector, which does not
scale with either the problem dimension or the network size, and is
not desirable from a privacy standpoint.~\cite{JL-CYT:16} also
exploits the positive definite property of the individual matrices and
requires the agents to know the state as well as the matrices of the
neighbors.~\cite{XW-SM:18} proposes a distributed algorithm without
any positive definiteness condition, but agents are allowed to
converge to different solutions.  Our approach here uses dynamic
average consensus~\cite{DPS-ROS-RMM:05b,SSK-BVS-JC-RAF-KML-SM:19} to
estimate certain non-distributed terms in a gradient-based algorithm
for the reformulated optimization problem.  We also draw inspiration
from~\cite{BT-BG:15,BT-BG:19} on distributed optimization to extend
our treatment to deal with unbalanced networks.  However, unlike the
aforementioned works where the desired solution is not an equilibrium
of the dynamics, requiring a diminishing time-varying stepsize-like
parameter to ensure convergence, here we make sure that any solution
of the linear equation is indeed an equilibrium of the proposed
dynamics. This enables us to employ Lyapunov stability analysis to
establish convergence and offers a framework to study robustness
against disturbances and errors. Our work~\cite{PS-JC:21-tcns}
requires bidirectional 2-hop communication. In contrast, the
distributed algorithms here require information exchange only with
immediate neighbors and work for arbitrary directed~graphs.

\footnotetext{We employ the following notation. $\real$, $\realpos$
  and $\integer$ denote the set of real numbers, positive real
  numbers, and integers, resp.  $|\mathcal{X}|$ denotes the
  cardinality of a set $\mathcal{X}$.  $\one$, $\zero$ and $I$ denote
  a vector or matrix of all ones and zeros, and an identity matrix of
  appropriate dimension, resp.  We let lowercase letters to denote
  vectors and uppercase letters to denote matrices.  $\|x\|$ and
  $\|A\|$ denote the 2-norm of a vector $x$ and the induced 2-norm of
  a matrix $A$, resp.  $\diag(x)$ denotes the diagonal matrix obtained
  after arranging the entries of the vector $x$ along the principal
  diagonal.  $A_{ij}$ denotes the $ij$th element of a matrix $A$,
  $A^\top$ its transpose, $A^{-1}$ its inverse (if it exists) and
  $\nulls(A)$ its null space.  $A \otimes B$ denotes the Kronecker
  product between two matrices $A$ and $B$.  Unless otherwise stated,
  $\xx \in \real^{mn}$ denotes the concatenated vector obtained after
  stacking the vectors $\{x_i\}_{i=1}^n \in \real^m$.  $A \succ \zero$
  and $A \succeq \zero$ imply that a matrix $A$ is positive definite
  and semidefinite, resp.  For a symmetric matrix $A$,
  $\lambda_{\max}(A)$ and $\lambda_{\min}(A)$ denote its maximum and
  minimum eigenvalue, resp.  Regardless of the multiplicity of
  eigenvalue 0, $\lambda_2(A)$ denotes the minimum non-zero eigenvalue
  of a positive semidefinite matrix $A$.  For two vectors $x,y \in
  \real^n$, $[x;y]$ denotes the concatenated vector containing the
  entries of $x$ and $y$, in that order, and $x > y$ means that the
  inequality holds elementwise.  }

\emph{Statement of Contributions:} We consider linear algebraic
equations where the coefficient matrices and constant vector for the
overall problem are given, respectively, by the summation of the
individual agents' coefficient matrices and constant vectors.  Our
starting point is the exact reformulation of this problem as a
constrained optimization problem.  Using the observation that the
optimal value of this optimization is zero, we reformulate it as
optimization of an unconstrained function, and propose a centralized
algorithm which works for weight-balanced networks and serves as a
reference for the design of distributed algorithms.  Using dynamic
average consensus, we then propose a distributed algorithm that does
not require the agent matrices to be positive definite, works for
time-varying weight-balanced networks and is guaranteed to converge to
a solution of the original problem exponentially fast.  Building on
the insights gained in establishing these results, we propose a
distributed algorithm that is not limited to weight-balanced networks
and is also guaranteed to converge to a solution of the linear
equation exponentially fast.

%\vspace*{-1ex}
\section{Preliminaries}\label{sec:prelim}
Here we review basic notions from graph
theory~\cite{FB-JC-SM:08cor,BT-BG:15,ZL-DZ:15} and dynamic average
consensus~\cite{DPS-ROS-RMM:05b,SSK-BVS-JC-RAF-KML-SM:19}.

\emph{Graph Theory:} Let $\G=(\V,\E,\A)$ be a weighted directed
graph (or digraph), with $\V$ as the set of vertices (or nodes) and
$\E\subseteq\V \times \V$ as the set of edges: $(v_i,v_j) \in \E$ iff
there is an edge from node $v_i$ to node $v_j$.  With $|\V|=n$, the
\emph{adjacency matrix} $\A \in \real^{n \times n}$ of $\G$ is such
that $\A_{ij}>0$ if $(v_i, v_j) \in \E$ and $\A_{ij}=0$, otherwise.  A
directed path is an ordered sequence of vertices such that any pair of
consecutive vertices is an edge. A digraph is strongly connected if
there is a directed path between any two distinct vertices. The out-
and in-degree of a node are, resp., the number of outgoing edges from
and incoming edges to it.  The weighted out-degree and weighted
in-degree of a node $v_i$ are $d^{\text{out}}(v_i)=\sum_{j=1}^{n}
\A_{ij}$ and $d^{\text{in}}(v_i)=\sum_{j=1}^{n} \A_{ji}$, resp.  The
\emph{out-degree matrix} $\D^{\text{out}} \in \real^{n \times n}$ and
 \emph{in-degree matrix} $\D^{\text{in}} \in \real^{n \times n}$
are diagonal matrices defined as $\D^{\text{out}}_{ii} =
d^{\text{out}} (v_i)$ and $\D^{\text{in}}_{ii}=d^{\text{in}}(v_i)$,
resp.  A graph is weight-balanced if $\D^{\text{out}}=\D^{\text{in}}$.
The \emph{Laplacian} $\Lap \in \real^{n \times n}$ is $\Lap =
\D^{\text{in}} - \A$.  All eigenvalues of $\Lap$ have nonnegative real
parts, $0$ is  simple with left eigenvector $\one$ iff
$\G$ is strongly connected, and $\Lap \one=\zero$ iff $\G$ is
weight-balanced iff $\Lap + \Lap^\top$ is positive semidefinite,
cf.~\cite[Theorem 1.37]{FB-JC-SM:08cor}.  If $\G$ is strongly
connected, it follows from~\cite[Lemma 3]{ZL-DZ:15} that there exists
a positive right eigenvector $\bar{v} \in \real^n$ associated to~$0$.
 
\emph{Dynamic Average Consensus:} Consider a group of $n \in~\integer_{> 1}$ 
agents communicating over a weight-balanced digraph
$\G$ with  Laplacian~$\Lap$. Each agent $i \in
\until{n}$ has a state $x_i \in \real$ and an input $z_i \in \real$.
The dynamic average consensus algorithm aims at making all the agents
track the average $\frac{1}{n} \sum_{i=1}^n z_i$ asymptotically. Here
we present the algorithm following~\cite{DPS-ROS-RMM:05b}, where it
was introduced for undirected graphs. Consider
\begin{align*}
  \dot{\xx}=-\Lap \xx + \dot{\zz}.
\end{align*}
If $\sum_{i=1}^n x_i(0)=\sum_{i=1}^n z_i(0)$ and the input $\zz$ is
bounded, then $x_i (t)\to \frac{1}{n} \sum_{i=1}^n z_i(t)$
as $t\rightarrow \infty$ for  $i \in \until{n}$, cf.~\cite{DPS-ROS-RMM:05b}.
%\vspace*{-0.5ex}

\section{Problem Formulation}\label{sec:prob}
Consider a group of $n$ agents interacting over a digraph that
seek to solve in a distributed way the linear algebraic equation
\begin{align}\label{eq:problem}
  \underbrace{\left( \sum\limits_{i=1}^n A_i \right)}_A x =
  \underbrace{\left(\sum\limits_{i=1}^n b_i\right)}_b,
\end{align}
where $x \in \real^m$ is the unknown solution vector, and $A_i \in
\real^{m \times m}$ and $b_i \in \real^m$ are the coefficient matrix
and constant vector corresponding to agent $i \in \until{n}$.  We
assume that~\eqref{eq:problem} has at least one solution.  The
formulation~\eqref{eq:problem} includes, as a particular case,
scenarios where each agent~$i$ knows only some rows of the coefficient
matrix~$A$ and constant vector~$b$.  Our approach consists of first
formulating~\eqref{eq:problem} as a system involving $n$ unknown
solution vectors, one per agent, and then reformulating it as a convex
optimization problem.  Based on this reformulation, we propose two
sets of (out-)distributed algorithms (where each agent only needs
information from its out-neighbors) to find the solutions
of~\eqref{eq:problem}.  We start by endowing each agent with its own
version $x_i \in \real^m$ of $x$. Then~\eqref{eq:problem} can be
equivalently written as
\begin{subequations}\label{eq:reform1}
  \begin{align}
    \sum\limits_{i=1}^n A_i x_i & = \sum\limits_{i=1}^n
    b_i, \label{eq:xcopy}
    \\
    x_i &= x_j \quad \forall i,j \label{eq:consensus}.
  \end{align}
\end{subequations}
Equation~\eqref{eq:consensus} ensures that $x_i=x$ for all the agents.
Clearly the set of equations~\eqref{eq:reform1} and the original
problem~\eqref{eq:problem} are equivalent.  Next we
formulate~\eqref{eq:reform1} as a convex optimization problem.
Consider the quadratic function $f: \real^{mn} \to \real$
\begin{align*}
  f(\xx) = \Big( \sum\limits_{i=1}^n (A_i x_i - b_i)\Big)^\top \Big(
  \sum\limits_{i=1}^n (A_i x_i - b_i) \Big) ,
\end{align*}
which is convex and attains its minimum over the solution set
of~$\eqref{eq:xcopy}$.  For convenience, we use $\LL=\Lap \otimes I$
and $f(\xx) = (\AAA \xx - \bb)^\top \onem \onem^\top (\AAA \xx-\bb)$,
where $\onem=\one \otimes I$, $\AAA \in \real^{mn \times mn}$ denotes
the block-diagonal matrix obtained after putting the matrices
$\{A_i\}_{i=1}^n$ along the principal diagonal, and $\bb=[b_1; \ldots;
b_n] \in \real^{mn}$.  If $\G$ is strongly connected, the solutions
of~\eqref{eq:reform1} are the same as the optimizers~of
\begin{equation}\label{eq:optimization}
  \begin{aligned}
    &\min\limits_{\xx} & & f(\xx)\\
    & \text{s.t.}  & & \LL^\top \xx = \zero.
  \end{aligned}
\end{equation}

\begin{remark}\longthmtitle{Distributed algorithmic solutions to
    optimization problem} 
  The problem~\eqref{eq:optimization} can be solved over an undirected
  graph by reformulating it using the techniques
  in~\cite{AC-JC:16-allerton} and employing the saddle-point dynamics,
  cf.~\cite{KA-LH-HU:58,AC-BG-JC:17-sicon}.  These dynamics involve
  terms of the form $\LL^\top$ and, to be implemented over a digraph,
  would need information from in- as well as out-neighbors and hence
  are not suitable for our setup.  Works that deal with distributed
  optimization under consensus constraints over digraphs, see
  e.g.~\cite{BG-JC:14-tac,BT-BG:15} and references therein, require
  the objective function to be separable, and therefore are not
  applicable either here.  \oprocend
\end{remark}

\section{Distributed Algorithms Over Weight-Balanced Networks}
We present distributed algorithms to solve~\eqref{eq:problem}
over weight-balanced networks.

\subsection{Centralized Algorithm}\label{sec:central}
We first introduce a centralized algorithm using the fact that the
objective function vanishes at the optimizers
of~\eqref{eq:optimization}.  Let
\begin{align}\label{eq:unconstrained}
  \min\limits_{\xx} \quad \frac{1}{2}\alpha \xx^\top (\LL + \LL^\top)
  \xx + \beta f(\xx),
\end{align}
where $\alpha, \beta >0$. Clearly,~\eqref{eq:optimization}
  and~\eqref{eq:unconstrained} have the same set of solutions if $\G$
  is strongly connected and weight-balanced. Since
problem~\eqref{eq:unconstrained} is unconstrained, one can use
gradient descent to find its optimizers. However, the gradient
$-\alpha (\LL+\LL^\top) \xx - \beta \AAA^\top \onem
  \onem^\top (\AAA \xx - \bb)$ of the objective function
in~\eqref{eq:unconstrained} involves terms with $\LL^\top$, whose
computation would require information from in-neighbors.  Instead, we
consider the following gradient-based dynamics
\begin{align}\label{eq:central}
  \dot{\xx} = -\alpha \LL \xx - \beta \AAA^\top \onem \onem^\top (\AAA
  \xx - \bb).
\end{align}
Whenever convenient, we refer to~\eqref{eq:central} as
$\grad$.  Note that the first term in the dynamics~\eqref{eq:central}
is distributed, meaning that each agent can implement it with
information from its out-neighbors.  The second term, however,
requires collective information from all the agents because of the
summation across the network.  Nevertheless, this algorithm  serves
as the basis for our distributed algorithm design in the next
section. 

The next result formally characterizes the equivalence
between the equilibria of~\eqref{eq:central} and the solutions
of~\eqref{eq:problem}.

\begin{lemma}\longthmtitle{Equivalence between~\eqref{eq:central}
    and~\eqref{eq:problem}}\label{lemma:equi}
  Let $\G$ be a strongly connected and weight-balanced digraph. Then
  for all $\alpha, \beta \in \realpos$, $\xx^*$ is an equilibrium
  of~\eqref{eq:central} if and only if $\xx^* = \one \otimes x^*$,
  where $x^* \in \real^m$ solves~\eqref{eq:problem}.
\end{lemma}
\begin{IEEEproof}
  The implication from right to left is immediate.  To prove the
  implication in the other direction, let $\bar{x} \in \real^m$ be a
  solution of~\eqref{eq:problem} and consider $\bar{\xx} = \one
  \otimes \bar{x}$.  Since $\xx^*$ and $\bar{\xx}$ are equilibria
  of~\eqref{eq:central},
  \begin{align}\label{eq:equilibrium}
    \alpha \LL (\xx^* - \bar{\xx}) + \beta \AAA^\top \onem \onem^\top
    \AAA (\xx^* - \bar{\xx}) = \zero.
  \end{align}
  Let $\QQ_{11}=\frac{1}{2}\alpha (\LL+\LL^\top)+\beta \AAA^\top \onem
  \onem^\top \AAA$.   Then~\eqref{eq:equilibrium} implies
  \begin{align*}
    (\xx^* - \bar{\xx})^\top \QQ_{11} (\xx^* - \bar{\xx}) =0.
  \end{align*}
  Since $\G$ is weight-balanced, $(\LL+\LL^\top) \succeq
    \zero$.  This along with the fact that $\AAA^\top \onem \onem^\top
    \AAA \succeq \zero$ implies $\LL^\top (\xx^* - \bar{\xx})=\zero$
  and $\onem^\top \AAA (\xx^* - \bar{\xx})=\zero$.  Therefore,
  $\xx^*=\one \otimes x^*$, for some $x^* \in~\real^m$ which satisfies
  $A x^* = A \bar{x} = b$, as claimed.
\end{IEEEproof}

The next result characterizes the convergence
of~\eqref{eq:central}.

\begin{proposition}\longthmtitle{Exponential stability
    of~\eqref{eq:central}}\label{prop:central}
  Let $\G$ be a strongly connected and weight-balanced digraph. Then
  for all $\alpha, \beta \in \realpos$, any trajectory
  of~\eqref{eq:central} converges exponentially to a point of the form
  $\xx^* = \one \otimes x^*$, where $x^* \in \real^m$
  solves~\eqref{eq:problem}.
\end{proposition}
\begin{IEEEproof}
  Consider a vector $\ww \in \real^{mn}$ in the null space of
  $\QQ_{11}$.  Using the same line of arguments as in the proof of
  Lemma~\ref{lemma:equi}, this
    implies that $\LL^\top \ww=\zero$ and $\onem^\top \AAA
  \ww=\zero$. Therefore, along~\eqref{eq:central},
  \begin{align*}
    \dot{\xx}^\top \ww= -(\alpha \xx^\top \LL^\top + \beta (\AAA \xx -
    \bb)^\top \onem \onem^\top \AAA) \ww=0.
  \end{align*}
  This means that the dynamics~\eqref{eq:central} are orthogonal to
  the null space of $\QQ_{11}$ and hence the component of $\xx$ in the
  null space of $\QQ_{11}$, say $\xx_{\nulls}$, remains constant.
  Given the initial condition $\xx(0)$, consider the particular
  equilibrium $\xx^*$ of~\eqref{eq:central} satisfying
  $\xx^*_{\nulls}=\xx(0)_{\nulls}$.  Since different equilibria differ
  only in their null space component,
   $\xx^*$ defined this way is unique.  Consider the Lyapunov function
  candidate $V : \real^{mn} \to \real$
  \begin{align*}
    V(\xx) = \frac{1}{2}(\xx - \xx^*)^\top (\xx - \xx^*).
  \end{align*}
  The Lie derivative of $V$ along the dynamics~\eqref{eq:central} is
  given by
  \begin{align*}
    \lie_{\grad}V=&-(\xx-\xx^*)^\top (\alpha \LL \xx + \beta \AAA^\top
    \onem \onem^\top (\AAA \xx - \bb))
    \\
    =&-(\xx-\xx^*)^\top \QQ_{11}(\xx-\xx^*) \le - 2\lambda_2(\QQ_{11})
    V.
  \end{align*}
  The last inequality follows from applying the
    Courant-Fischer theorem~\cite[Theorem 4.2.11]{RAH-CRJ:85} together
    with the fact that $(\xx-\xx^*)^\top \ww=0$ as $\xx_{\nulls}$ is
    constant.  Using the monotonicity theorem~\cite[Corollary
    4.3.3]{RAH-CRJ:85}, we further have
  \begin{align*}
    \lie_{\grad} V \le - 2 \min\left\{\frac{1}{2} \alpha
      \lambda_2(\Lap+\Lap^\top), \beta \lambda_2(\AAA^\top \onem
      \onem^\top \AAA)\right\}V.
  \end{align*}
  Hence, the dynamics~\eqref{eq:central} is exponentially stable with
  a rate depending on $\alpha, \beta, \Lap$ and $\{A_i\}_{i=1}^n$.
\end{IEEEproof}
%\vspace*{-2ex}
\subsection{Distributed Algorithm}\label{sec:distributed}
We present a distributed algorithm to find a solution
of~\eqref{eq:problem}, which is based on the centralized
algorithm~\eqref{eq:central} and involves employing dynamic
  average consensus (cf.~Section~\ref{sec:prelim}) to estimate the aggregate
  $\onem^\top (\AAA \xx - \bb)$.  Formally,
\begin{subequations}\label{eq:distributed}
  \begin{align}
    \dot{\xx}=&-\alpha \LL \xx - n \beta \AAA^\top
    \yy, \label{eq:distributed_a}  
    \\
    \dot{\yy}=&-\alpha \AAA \LL \xx - n \beta \AAA \AAA^\top \yy -
    \gamma \LL \yy, \label{eq:distributed_b}
  \end{align}
\end{subequations}
with design parameter $\gamma > 0$. Here, each agent $i \in
  \until{n}$ updates  $y_i \in \real^m$ which estimates the
  average mismatch $\frac{1}{n} \onem^\top (\AAA \xx - \bb)$. The
dynamics~\eqref{eq:distributed} is distributed as each agent just
needs to know its state and that of its out-neighbors. Whenever
convenient, we refer to it as $\gdac$.  The following result
characterizes the equilibria of~\eqref{eq:distributed} and shows that
the total deviation from the average mismatch is conserved.

\begin{lemma}\longthmtitle{Equilibria of~\eqref{eq:distributed}
    and invariance of total deviation}\label{lemma:equi_distributed}
  Let $\G$ be a strongly connected and weight-balanced digraph.
  Then, if $(\xx^*,\zero)$ is an equilibrium of~\eqref{eq:distributed}
  then $\xx^* = \one \otimes x^*$, where $x^* \in \real^m$.  Moreover,
  for all $\alpha, \beta , \gamma \in \realpos$, $\onem^\top (\yy -
  \AAA \xx) $ remains constant along the evolution
  of~\eqref{eq:distributed}.
\end{lemma}
\begin{IEEEproof}
  Let $(\xx^*,\zero)$ be an equilibrium
  of~\eqref{eq:distributed}. From~\eqref{eq:distributed_a}, it follows
  that $\LL \xx = \zero$, and hence $\xx^*=\one \otimes x^*$ for some $x^*
  \in \real^m$, establishing the first statement. Now, consider the
  derivative $\onem^\top (\dot{\yy} - \AAA \dot{\xx})= - \gamma
  \onem^\top \LL \yy = \zero$. Hence, $\onem^\top (\yy - \AAA \xx)$ is
  conserved along the evolution of~\eqref{eq:distributed}.
\end{IEEEproof}

\begin{remark}\longthmtitle{Distributed initialization of the $\gdac$
    algorithm}\label{remark:init}
  From Lemma~\ref{lemma:equi_distributed}, we observe that in order
  for a trajectory of~\eqref{eq:distributed} to converge to an
  equilibrium of the form $(\xx^*, \yy^*) = (\one \otimes x^*, \zero)
  $, where $x^* \in \real^m$ solves~\eqref{eq:problem}, its initial
  condition must satisfy
     $\onem^\top \yy(0) = \onem^\top (\AAA \xx(0) - \bb).$
   This could be implemented in a distributed way if each agent $i \in
  \until{n}$ chooses its initial states satisfying $y_i(0)=A_i x_i(0)
  - b_i$. One trivial selection, for example, is $\xx(0)=\zero$ and
  $\yy(0)=-\bb$.  \oprocend
\end{remark}

The next result characterizes the convergence
of~\eqref{eq:distributed}.

\begin{theorem}\longthmtitle{Exponential stability
    of~\eqref{eq:distributed} over balanced
    networks}\label{thm:distributed} 
  Let $\G$ be a strongly connected and weight-balanced digraph
  and assume 
$\nulls(A)   \subseteq \nulls(A_i)$, for all $i \in
  \until{n}$. Let $\alpha, \beta \in \realpos $ and define 
  \begin{align*}
  \bar{\gamma} =   \max\left\{ \frac{2}{\lambda_2{(\Lap+\Lap^\top)}}
    \lambda_{\max}\left(\frac{\QQ_{12}^\top
      \QQ_{12}}{\lambda_2(\QQ_{11})}- n \beta \AAA \AAA^\top\right), 0 \right\},
  \end{align*}
  where $\QQ_{11} =\frac{1}{2}\alpha (\LL+\LL^\top)+\beta \AAA^\top \onem
  \onem^\top \AAA$ and
  $\QQ_{12}=\frac{1}{2} (n \beta \AAA^\top + \alpha \LL^\top \AAA^\top
  + \beta \AAA^\top \onem \onem^\top \AAA \AAA^\top)$.  Then, for all
  $\gamma \in (\bar{\gamma}, \infty)$, any trajectory
  of~\eqref{eq:distributed} with initial condition satisfying
  $\onem^\top \yy(0)=\onem^\top (\AAA \xx(0)- \bb)$ converges
  exponentially to $(\xx^*, \zero)$, where $\xx^*=\one \otimes x^*$
  and $x^* \in \real^m$ solves~\eqref{eq:problem}.
\end{theorem}
\begin{IEEEproof}
    Define the error variable
  \begin{align}\label{eq:e}
    \ee=\yy-\frac{1}{n}\onem \onem^\top (\AAA \xx - \bb),
  \end{align}
  measuring the difference between the agents' estimates and the
  actual value of average mismatch. Note that
  \begin{align*}
    \dot{\ee} = -\alpha \PiPi \AAA \LL \xx - n \beta \PiPi \AAA \AAA^\top \yy -
    \gamma \LL \yy,
  \end{align*}
  where $\PiPi= I - \frac{1}{n} \onem \onem^\top$.  Rewriting
  \eqref{eq:distributed} in terms of $\xx$ and $\ee$,
  \begin{subequations}\label{eq:transformed}
    \begin{align}\label{eq:transformed_a}
      \dot{\xx}&=-\alpha \LL \xx - \beta \AAA^\top \onem \onem^\top
      (\AAA \xx - \bb) - n \beta \AAA^\top \ee,
      \\
      \dot{\ee}&=-\alpha \PiPi \AAA \LL \xx - \beta \PiPi \AAA
      \AAA^\top \onem \onem^\top (\AAA \xx - \bb)
      \\
      & \quad - n \beta \PiPi \AAA \AAA^\top \ee -\gamma \LL
      \ee. \notag
    \end{align}
  \end{subequations}
  From the proof of Proposition~\ref{prop:central}, we know that if
  $\ww \in \real^{mn}$ is in the null space of $\QQ_{11}$, then
  $\LL^\top \ww=\zero$ and $\onem^\top\AAA \ww=~\zero$. Therefore, $\ww
  = \one \otimes w$, where $w \in \real^m$ belongs to $w \in
  \nulls(A)$. By hypothesis, $A_i w = \zero$ for all $i \in
  \until{n}$.  Therefore, from~\eqref{eq:transformed_a},
  $\dot{\xx}^\top \ww =0$, and the $\xx$ component of the
  equilibrium $(\xx^*,\yy^*)$ of~\eqref{eq:distributed} satisfies
  $\xx^*_{\nulls}=\xx(0)_{\nulls}$ and is unique.  With the
  initialization of the statement, it follows from
  Lemma~\ref{lemma:equi_distributed} that $\yy^*=\one \otimes
  \frac{1}{n} \onem^\top (\AAA \xx^* - \bb)$.  Substituting this value
  of $\yy^*$ in~\eqref{eq:distributed_a} and following the proof of
  Lemma~\ref{lemma:equi}, one can establish that the corresponding
  equilibrium is of the form $(\one \otimes x^*, \zero)$, where $x^*
  \in \real^m$ is a solution of~\eqref{eq:problem}.  Consider the
  Lyapunov function candidate $V_2 : \real^{2mn} \to \real$
  \begin{align*}
    V_2(\xx,\ee)=\frac{1}{2}(\xx-\xx^*)^\top(\xx-\xx^*)+\frac{1}{2}\ee^\top
    \ee.
  \end{align*}
  The Lie derivative of $V_2$ along $\eqref{eq:transformed}$ is given
  by
  \begin{align*}
    \lie_{\gdac}V_2=&-(\xx-\xx^*)^\top(\alpha \LL \xx+ \beta \AAA^\top
    \onem \onem^\top (\AAA \xx - \bb)) \\ &- n \beta (\xx -
    \xx^*)^\top \hspace*{-0.5ex} \AAA^\top \hspace*{-0.5ex} \ee
    -\ee^\top \PiPi \AAA (\alpha  \LL \xx + n \beta \AAA^\top \hspace*{-0.5ex} \ee)\\
    & -\ee^\top (\beta \PiPi \AAA \AAA^\top \onem \onem^\top (\AAA \xx
    - \bb)+\gamma \LL \ee) \\
    =&-\begin{bmatrix} \xx-\xx^* \\
      \ee \end{bmatrix}^\top \begin{bmatrix} \QQ_{11} & \QQ_{12} \\
      \QQ_{12}^\top & \QQ_{22}\end{bmatrix} \begin{bmatrix} \xx-\xx^*
      \\ \ee \end{bmatrix},
  \end{align*} 
  where $\QQ_{22}=\frac{1}{2}\gamma (\LL+\LL^\top) + n \beta \AAA
  \AAA^\top$ and we have used the fact that due to the mentioned
  initialization, $\onem^\top \ee=\zero$ from
  Lemma~\ref{lemma:equi_distributed}.  Since $\xx_{\nulls}$
    is constant, $(\xx-\xx^*)^\top \ww=0$ and from the
    Courant-Fischer theorem~\cite[Theorem 4.2.11]{RAH-CRJ:85},
  \begin{align*}
    -(\xx-\xx^*)^\top \QQ_{11} (\xx-\xx^*) \leq - \lambda_2(\QQ_{11})(\xx-\xx^*)^\top
     (\xx-\xx^*).
  \end{align*}
  Also, since $\onem^\top \ee=\zero$ and $\G$ is
    weight-balanced, it again follows from the Courant-Fischer theorem
    that
  \begin{align*}
    -\ee^\top \QQ_{22} \ee \leq -\frac{1}{2} \gamma 
    \lambda_2(\Lap + \Lap^\top)\ee^\top \ee - n \beta \ee^\top \AAA \AAA^\top
    \ee .
  \end{align*}
  Therefore, we can upper bound the Lie derivative as
  \begin{align*}
    \lie_{\gdac} V_2 \le -
    \begin{bmatrix}
      \xx-\xx^* \\
      \ee
    \end{bmatrix}^\top \underbrace{\begin{bmatrix}
        \lambda_2(\QQ_{11})I & \QQ_{12}
        \\
        \QQ_{12}^\top &
        \bar{\QQ}_{22}\end{bmatrix}}_{\bar{\QQ}} \begin{bmatrix}
      \xx-\xx^* \\ \ee \end{bmatrix},
  \end{align*}
  where $\bar{\QQ}_{22}= \frac{1}{2} \gamma \lambda_2(\Lap +
  \Lap^\top)I + n \beta \AAA \AAA^\top$.  Next, we examine the
  positive definiteness of $\bar{\QQ}$.  Using the Schur
  complement~\cite{SB-LV:09}, $\bar{\QQ} \succ \zero$ iff
  \begin{align*}
    \frac{1}{2} \gamma \lambda_2(\Lap+\Lap^\top)I + n \beta \AAA
    \AAA^\top - \frac{1}{ \lambda_2(\QQ_{11})}\QQ_{12}^\top
    \QQ_{12}\succ \zero.
  \end{align*}
  Hence, $\bar{\QQ} \succ \zero$ if $\gamma > \bar{\gamma}$, and
  $\lie_{\gdac}V_2 \leq -2\lambda_{\min}(\bar{\QQ})V_2$.
\end{IEEEproof}

The null space condition in Theorem~\ref{thm:distributed} makes sure
that $\xx^*_{\nulls}$ remains invariant along the evolution
of~\eqref{eq:distributed} and all the agents approach the solution
of~\eqref{eq:problem} closest to $\xx(0)$.  This condition is
automatically satisfied if the matrix $A$ is full rank, or in other
words, equation~\eqref{eq:problem} has a unique solution.  We believe
(and simulations also suggest) that if this condition is not
satisfied, the $\xx$ component of the dynamics still converges to a
solution of~\eqref{eq:problem}.

\begin{remark}\longthmtitle{Lower bound on $\gamma$} 
  The lower bound $\bar{\gamma}$ in Theorem~\ref{thm:distributed} is
  conservative in general.
    In fact, the algorithm may converge even if this condition is not
  satisfied, something that we have observed in simulation.  Note also
  that although $\alpha$ and $\beta$ are free parameters, they should
  still be carefully chosen as $\bar{\gamma}$ depends on them.
  \oprocend
\end{remark}

The result above can be extended to time-varying networks.  In case
$\G(t)$ is time-varying, the algorithm in~\eqref{eq:distributed} reads
as
\begin{subequations}\label{eq:time_varying}
  \begin{align}
    \dot{\xx}&=-\alpha \LL(t) \xx - n \beta \AAA^\top \yy,
    \\
    \dot{\yy}&=-\alpha \AAA \LL(t) \xx - n \beta \AAA \AAA^\top \yy -
    \gamma \LL(t) \yy.
  \end{align}
\end{subequations}
The next result formally characterizes the convergence
of~\eqref{eq:time_varying}. Its proof is similar to that of
Theorem~\ref{thm:distributed} and hence omitted.

\begin{theorem}\longthmtitle{Exponential stability
    of~\eqref{eq:time_varying} over time-varying balanced
    networks}\label{thm:time_varying}
  Let $\{\G(t)\}_{t=0}^{\infty}$ be a sequence of strongly connected
  and weight-balanced digraphs with uniformly bounded edge
  weights (i.e., there exists $a \in (0,\infty)$ such that
    $\A_{ij}(t)<a$ for all $(i,j)$ and $t \ge 0$), and assume
  $\nulls(A) \subseteq \nulls(A_i)$, for all $i \in
  \until{n}$. Let $\alpha, \beta \in \realpos $ and define
  $\bar{\gamma}(t)$ as
  \begin{align*}
    \max \hspace*{-0.5ex} \left\{ \hspace*{-0.5ex}
      \frac{2}{\lambda_2(\Lap(t) \hspace*{-0.5ex} + \hspace*{-0.5ex}
        \Lap(t)^\top)}    
    \lambda_{\max} \hspace*{-0.5ex}  \left(\hspace*{-0.1ex} \frac{\QQ_{12}(t)^\top \hspace*{-0.4ex}
      \QQ_{12}(t)}{\lambda_2(\QQ_{11}(t))} \hspace*{-0.5ex} -  \hspace*{-0.5ex} n \beta \AAA \AAA^\top \hspace*{-0.5ex}  \right) \hspace*{-0.6ex}, 0 \hspace*{-0.5ex} \right\} \hspace*{-0.6ex},
  \end{align*}
  where $\QQ_{11}(t)=\frac{1}{2}\alpha (\LL(t)+ \LL(t)^\top)+ \beta
  \AAA^\top \onem \onem^\top \AAA$ and $\QQ_{12}(t)=\frac{1}{2} (n
  \beta \AAA^\top + \alpha \LL(t)^\top \AAA^\top + \beta \AAA^\top
  \onem \onem^\top \AAA \AAA^\top)$.  Then for all $\gamma \in
  (\hat{\gamma}, \infty)$, where $\hat{\gamma} = \sup\limits_{t \ge 0}
  \bar{\gamma} (t)$, any trajectory of~\eqref{eq:time_varying} with
  initial conditions $\onem^\top \yy(0)=\onem^\top (\AAA \xx(0)- \bb)$
  converges exponentially to $(\xx^*, \zero)$, where $\xx^*=\one
  \otimes x^*$ and $x^* \in \real^m$ solves~\eqref{eq:problem}.
\end{theorem}

\section{Distributed Algorithm Over Unbalanced
  Networks}\label{sec:unbalanced}
In this section, we extend our approach to solve
problem~\eqref{eq:problem} over graphs that are not necessarily
balanced.  In those scenarios, since $\Lap\one \ne \zero$, the
one-to-one correspondence between the desired equilibria
of~\eqref{eq:central} or~\eqref{eq:distributed} and the solutions
of~\eqref{eq:problem} does not hold anymore.  To overcome this, we
propose
\begin{subequations}\label{eq:distributed_unbalanced}
  \begin{align}
    \dot{\xx}&=-\alpha \LL \bar{\VV} \xx - n \beta \AAA^\top
    \yy, \label{eq:distributed_unbalanced_a}
    \\
    \dot{\yy}&=-\alpha \AAA \LL \bar{\VV} \xx - n \beta \AAA \AAA^\top
    \yy - \gamma \LL \bar{\VV}
    \yy, \label{eq:distributed_unbalanced_b}
  \end{align}
\end{subequations}
where $\bar{\VV} = \diag(\bar{\vv})$, $\bar{\vv}=\one \otimes \bar{v}$, and $\bar{v}$ is a positive right
eigenvector with eigenvalue $0$ of $\Lap$.
Exponential stability
of~\eqref{eq:distributed_unbalanced} can be established by
interpreting $\Lap \cdot \diag(\bar{v})$
as the Laplacian of a weight-balanced graph and then following the
same steps as in the proof of Theorem~\ref{thm:distributed}, but we
omit it here for reasons of space.
Although~\eqref{eq:distributed_unbalanced} is distributed, it assumes
that agents have a priori knowledge of the corresponding entries of
$\bar{v}$ which might be limiting in practice.  To deal with this limitation, 
we propose an  algorithm that does not require
such knowledge by augmenting~\eqref{eq:distributed_unbalanced} with an
additional dynamics converging to $\bar{\vv}$,
\begin{subequations}\label{eq:distributed_v}
  \begin{align}
    \dot{\xx}&=-\alpha \LL \VV \xx - n \beta \AAA^\top
    \yy, \label{eq:distributed_v_a}
    \\
    \dot{\yy}&=-\alpha \AAA \LL \VV \xx -  n \beta \AAA \AAA^\top \yy
    - \gamma \LL \VV \yy, \label{eq:distributed_v_b} 
    \\
    \dot{\vv}&=-\LL \vv, \label{eq:distributed_v_c}
  \end{align}
\end{subequations}
where $\VV=\diag(\vv)$.  Whenever convenient, we refer to
dynamics~\eqref{eq:distributed_v} as $\dist$. Note that, unlike all
the dynamics discussed so far, $\dist$ is nonlinear.

\begin{remark}\longthmtitle{Distributed nature of~\eqref{eq:distributed_v}}
    The dynamics~\eqref{eq:distributed_v} is out-distributed, but
  requires each agent $i \in \until{n}$ to have knowledge of its
  in-degree because $\Lap = \D^{\textnormal{in}}-\A$ and the graph is
  not weight-balanced. If we use instead the out-Laplacian
  $\Lap=\D^{\textnormal{out}}-\A$, then one could still define an
  equivalent algorithm for~\eqref{eq:distributed_unbalanced} with $\LL
  \bar{\VV}$ replaced by $\bar{\VV} \LL$,
  but~\eqref{eq:distributed_v_c} would look like $\dot{\vv}=-\LL^\top
  \vv$, which would require state information from in-neighbors too.
  \oprocend
\end{remark}

The next result characterizes the convergence
of~\eqref{eq:distributed_v}.

\begin{theorem}\longthmtitle{Exponential stability of~\eqref{eq:distributed_v}
    over unbalanced networks}\label{thm:distributed_v}
  Let $\G$ be a strongly connected digraph and assume
  $\nulls(A) \subseteq \nulls(A_i)$, for all $i \in
  \until{n}$. Let $\alpha, \beta \in \realpos $ and define
  \begin{align*}
    \bar{\gamma} \hspace*{-0.5ex}= \hspace*{-0.5ex} \max
    \hspace*{-0.5ex} \left\{ \frac{2}{\lambda_2{(\Lap \bar{V} +
          \bar{V} \Lap^\top)}} \lambda_{\max} \hspace*{-0.5ex}\left(
        \hspace*{-0.5ex} \frac{\QQ_{12}^\top
          \QQ_{12}}{\lambda_2(\QQ_{11})}- n \beta \AAA \AAA^\top
        \hspace*{-0.5ex} \right) \hspace*{-0.7ex}, 0 \right\}
    \hspace*{-0.5ex},
  \end{align*}
  where $\QQ_{11}=\frac{1}{2}(\alpha \LL \bar{\VV} + \bar{\VV}
  \LL^\top) + \beta \AAA^\top \onem \onem^\top \AAA$,
  $\QQ_{12}=\frac{1}{2} (n \beta \AAA^\top + \alpha \bar{\VV} \LL^\top
  \AAA^\top + \beta \AAA^\top \onem \onem^\top \AAA \AAA^\top)$,
  $\bar{v}$ is the positive eigenvector with eigenvalue $0$ of $\Lap$
  satisfying $\one^\top \bar{v}=1$, and $\bar{V}=\diag(\bar{v})$.
  Then, for all $\gamma \in (\bar{\gamma}, \infty)$, any trajectory
  of~\eqref{eq:distributed_v} with initial condition satisfying
  $\onem^\top \yy(0)=\onem^\top (\AAA \xx(0)- \bb)$ and $\vv(0) =
  \frac{1}{n} \one$, converges exponentially to $(\xx^*, \zero,
  \bar{\vv})$, where $\xx^*=\one \otimes x^*$ and $x^* \in \real^m$
  solves~\eqref{eq:problem}, and $\bar{\vv}=\one \otimes \bar{v}$.
\end{theorem}
\begin{IEEEproof}
    From~\cite[Proposition 2.2]{BT-BG:19}, we have that $\vv(t) > \zero$
  for all $t \ge 0$. Also, since $\one^\top \Lap =\zero$, $\one^\top
  \vv$ is conserved along the evolution
  of~\eqref{eq:distributed_v_c}. Hence $\vv(t) \to \bar{\vv}$
  exponentially fast with a rate determined by the non-zero eigenvalue
  of $\Lap$ with the smallest real part.
  Let us interpret the
  dynamics~\eqref{eq:distributed_v_a}-\eqref{eq:distributed_v_b} as
  the dynamics~\eqref{eq:distributed_unbalanced} with some disturbance
  $\dd(t)$ defined by
  \begin{align*}
    \dd= \begin{bmatrix} \dd^{\xx} \\
      \dd^{\yy} \end{bmatrix}=   \begin{bmatrix} -\alpha \LL (\VV -
      \bar{\VV}) \xx  \\ 
      - \alpha \AAA \LL (\VV - \bar{\VV}) \xx - \gamma \LL (\VV -
      \bar{\VV} ) \yy \end{bmatrix},
  \end{align*}
  which goes to $\zero$ as $t \to \infty$.  Consider a vector $\ww \in
  \nulls(\QQ_{11})$. Then as in the proof of
  Theorem~\ref{thm:distributed}, $\ww = \one \otimes w$, where $w \in
  \nulls(A)$ and by hypothesis, $A_i w = \zero$ for all $i
  \in \until{n}$. Since $\onem^\top \LL=\zero$, therefore,
    $\ww^\top \dd^{\xx} = 0$ and we still have $\ww^\top
  \dot{\xx}=0$, and the $\xx$ component of the equilibrium
    $(\xx^*, \yy^*, \bar{\vv})$ of~\eqref{eq:distributed_v} satisfies
    $\xx^*_{\nulls}= \xx(0)_{\nulls}$ and is unique.  With the
    initialization of the statement and following the same steps as in
    the proof of Lemma~\ref{lemma:equi_distributed}, one can establish
    that $\yy^* = \one \otimes \frac{1}{n}\onem^\top (\AAA \xx^* -
    \bb)$.  Substituting this value of $\yy^*$
    in~\eqref{eq:distributed_v_a} and following the proof of
    Lemma~\ref{lemma:equi}, one can establish that the corresponding
    equilibrium is of the form $(\one \otimes x^*, \zero, \bar{\vv})$,
    where $x^* \in \real^m$ is a solution of~\eqref{eq:problem}.
  Consider now the Lyapunov function candidate $V_3 : \real^{3mn} \to
  \real$
  \begin{align*}
    V_3(\xx,\ee,\vv)=V_2(\xx,\ee)+ \frac{\delta}{2}(\vv -
    \bar{\vv})^\top \PP (\vv - \bar{\vv}),
  \end{align*}
  where $\delta > 0$, $\PP=\bar{\VV}^{-1}$, $\ee$ is defined as in~\eqref{eq:e}, and $V_2$ is the same
  function as in the proof of Theorem~\ref{thm:distributed}.
The Lie derivative of $V_3$ along
  $\eqref{eq:distributed_v}$ is given by
  \begin{align*}
    \lie_{\dist}V_3 \hspace*{-0.2ex} = & \hspace*{-0.2ex} -\begin{bmatrix} \xx-\xx^* \\
      \ee \end{bmatrix}^\top \hspace*{-1ex} \begin{bmatrix} \QQ_{11} & \QQ_{12} \\
      \QQ_{12}^\top & \QQ_{22}\end{bmatrix} \hspace*{-0.5ex}  \begin{bmatrix} \xx-\xx^*
      \\ \ee \end{bmatrix}  \hspace*{-0.5ex} + \hspace*{-0.5ex} (\xx-\xx^*)^\top \hspace*{-0.5ex}
    \dd^{\xx} \\
& + \ee^\top \dd^{\ee}   - \delta (\vv - \bar{\vv})^\top (\LL^\top \PP+
    \PP \LL) (\vv - \bar{\vv}) ,
  \end{align*}
  where $\dd^{\ee}= -\alpha \PiPi \AAA \LL (\VV- \bar{\VV}) \xx -
  \gamma \LL (\VV - \bar{\VV}) \ee$, and $\QQ_{22}=\frac{1}{2}\gamma
  (\LL \bar{\VV}+ \bar{\VV} \LL^\top) + n \beta \AAA \AAA^\top$.
  Interestingly, $\Lap \bar{V}$ can be interpreted as the
    Laplacian of a weight-balanced graph and as a result, $\LL
    \bar{\VV} + \bar{\VV} \LL^\top \succeq \zero$ implying that
    $\LL^\top \PP + \PP \LL \succeq \zero$.  Once again, following
  Lemma~\ref{lemma:equi_distributed}, one can establish that with the
  initialization of the statement, $\onem^\top \ee=\zero$ and
  therefore using the Courant-Fischer theorem~\cite[Theorem
    4.2.11]{RAH-CRJ:85} together with the fact that $(\xx-\xx^*)^\top
    \ww=0$ due to invariance of $\xx_{\nulls}$, we can upper bound
  the Lie derivative as
  \begin{align*}
    \lie_{\dist} V_3 \hspace*{-0.5ex}   \le&
    -\begin{bmatrix} \xx-\xx^* \\ \ee\end{bmatrix}^\top \hspace*{-1ex}
    \underbrace{\begin{bmatrix} \lambda_2(\QQ_{11})I \hspace*{-0.5ex} & \QQ_{12}\\ 
        \QQ_{12}^\top &
        \bar{\QQ}_{22}\end{bmatrix}}_{\bar{\QQ}} \hspace*{-0.5ex} \begin{bmatrix}
      \xx-\xx^* \\ \ee \end{bmatrix}     \\
    & \hspace*{-2ex} +\hspace*{-0.2ex} \alpha \|\xx-\xx^*\| \| \LL\| \|\vv - \bar{\vv}\| (\|\xx - \xx^*\|+ \|\xx^*\|) \\
&  \hspace*{-2ex}  +\hspace*{-0.2ex} \alpha \|\ee\| \| \PiPi \AAA \LL \| \|\vv - \bar{\vv}\|(\|\xx - \xx^*\|+ \|\xx^*\|) \\
 & \hspace*{-2ex}  +\hspace*{-0.2ex} \gamma \|\ee\| \|\LL\| \|\vv - \bar{\vv}\| \|\ee\|\hspace*{-0.5ex}  - \hspace*{-0.5ex} \delta \lambda_2(\LL^\top \hspace*{-0.5ex} \PP +
    \PP \LL)\|\vv - \bar{\vv}\|^2 \hspace*{-1ex}  ,
  \end{align*}
  where $\bar{\QQ}_{22}= \frac{1}{2} \gamma \lambda_2(\Lap \bar{V} +
  \bar{V} \Lap^\top)I + n \beta \AAA \AAA^\top$.  Define
  $\zz=[\|\xx-\xx^*\|;\|\ee\|;\|\vv-\bar{\vv}\|]$.  If $\gamma >
  \bar{\gamma}$, then $\bar{\QQ} \succ \zero$ and from the
    Courant-Fischer theorem, we have
  \begin{align*}
    \lie_{\dist} V_3 \le - \zz^\top \hspace*{-1ex} \underbrace{\begin{bmatrix} \lambda_{\min} (\bar{\QQ}) & 0 & \hat{\QQ}_{13}(\zz)\\
        0 & \lambda_{\min} (\bar{\QQ}) & \hat{\QQ}_{23}(\zz) \\
        \hat{\QQ}_{13}(\zz) & \hat{\QQ}_{23}(\zz) & \delta \lambda_2(\LL^\top \PP + \PP \LL)\end{bmatrix}}_{\hat{\QQ}(\zz)} \hspace*{-0.7ex}  \zz,
  \end{align*}
  where $\hat{\QQ}_{23}(\zz) = -\frac{1}{2} \alpha \|\PiPi \AAA \LL\|
  (\zz+ \|\xx^*\|) - \frac{1}{2}\gamma \|\LL\| \zz$ and
  $\hat{\QQ}_{13}(\zz)= -\frac{1}{2}\alpha \|\LL\| (\zz+ \|\xx^*\|)$.
  Using the Schur complement, one can verify that for a given value of
  $\zz$, $\hat{\QQ}(\zz) \succ \zero$ iff $\delta > \bar{\delta}(\zz)=
  \dfrac{1}{\lambda_{\min}(\bar{\QQ}) \lambda_2 (\LL^\top \PP + \PP
    \LL)} (\hat{\QQ}_{13}(\zz)^2 + \hat{\QQ}_{23}(\zz)^2).$ Hence, if
  $\delta > \bar{\delta}(\zz(0))$, then $\lie_{\dist} V_3 \leq -
  \lambda_{\min} (\hat{\QQ}(\zz(0)))\zz^\top \zz$.  This along with
  the fact that $\frac{1}{2} \min\{1, \delta \lambda_{\min}(\PP)\}
  \|\zz\|^2 \le V_3 \le \frac{1}{2} \max\{1, \delta
  \lambda_{\max}(\PP)\} \|\zz\|^2$, implies that $V_3$ satisfies the
  hypotheses of~\cite[Theorem 4.10]{HKK:02} for exponential stability.
\end{IEEEproof}

The exponential convergence of algorithms~\eqref{eq:central}
  and~\eqref{eq:distributed} for weight-balanced graphs,
  and~\eqref{eq:distributed_unbalanced} for unbalanced graphs follows
  from their linear nature.  For algorithm~\eqref{eq:distributed_v},
  exponential convergence could be attributed to the fact that the
  dynamics~\eqref{eq:distributed_v_c} converge exponentially and
  hence, after some
  time,~\eqref{eq:distributed_v_a}-\eqref{eq:distributed_v_b}
  and~\eqref{eq:distributed_unbalanced} are essentially the same.  

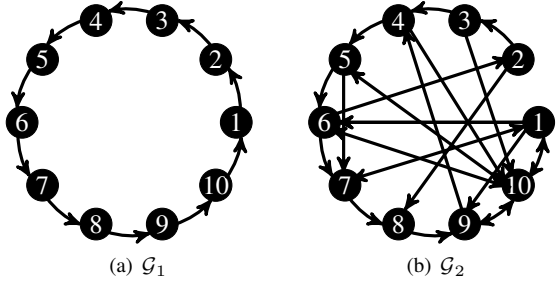
\begin{figure}[htb]
  \centering
  \tikzset{filled node/.style={draw, fill, minimum width=12pt, inner
      sep=0pt, circle}}
  \subfigure[$\G_1$]
  {\centering
    \begin{tikzpicture}[scale=0.9,every node/.style=filled node]
      \def \n {10}
      \def \radius {45pt}
      \foreach \s in {1,...,\n}
      {
        \node at ({360/\n * (\s - 1)}:\radius) {\textcolor{white}{\s}};
        \draw[->,>=stealth', very thick] ({360/\n * (\s - 1)+6}:\radius+2) 
        arc  ({360/\n * (\s - 1)}:{360/\n * (\s)-4}:\radius-11);
        }
    \end{tikzpicture} \label{fig:ring}
  }
\hspace*{7pt}
\subfigure[$\G_2$]
    {\centering
    \begin{tikzpicture}[scale=0.9,every node/.style=filled node]
      \def \n {10}
      \def \radius {45pt}
     \foreach \s in {1,...,\n}
{
\node at ({360/\n * (\s - 1)}:\radius) {\textcolor{white}{\s}};
}
\foreach \s in {2,...,8}
      {\draw[->,>=stealth', very thick] ({360/\n * (\s - 1)+6}:\radius+2) 
        arc ({360/\n * (\s - 1)}:{360/\n * (\s)-4}:\radius-11);
}
\draw[<->,>=stealth', very thick] ({360/\n * (\n - 2)+6}:\radius+2) arc ({360/\n * (\n - 2)}:{360/\n * (\n-1)-4}:\radius-11);
 \draw[<->,>=stealth', very thick] ({360/\n * (\n - 1)+6}:\radius+2) arc ({360/\n * (\n - 1)}:{360/\n * (\n)-4}:\radius-11);
     \draw[->,>=stealth', very thick] (1.4,0) -- (-1.4,0);
  \draw[<->,>=stealth', very thick] (1.4,-0.05) -- (-1.2,-0.85);
   \draw[->,>=stealth', very thick] (1.4,-0.1) -- (0.5,-1.34);
   \draw[->,>=stealth', very thick] (1.4,0) -- (-1.4,0);
 \draw[->, >=stealth', very thick] (1.15,0.85) -- (-0.4,-1.3);
\draw[->, >=stealth', very thick] (0.5,1.3) -- (1.2,-0.8);
 \draw[->, >=stealth', very thick] (-0.4,1.5) -- (1.1,-0.9);
  \draw[->, >=stealth', very thick] (-1.3,0.75) -- (-1.3,-0.75);
 \draw[<->, >=stealth', very thick] (-1.5,-0.1) -- (1.1,-0.95);
 \draw[->, >=stealth', very thick] (0.45,-1.3) -- (-0.4,1.4);
\draw[->, >=stealth', very thick] (-1.5,0.1) -- (1.15,0.95);
 \draw[->, >=stealth', very thick] (1.17,-0.97) -- (-1.23,0.78);
\end{tikzpicture}
\label{fig:ring+edges}
  }
  \caption{Communication topologies among the agents.  The edge
    weights are adjusted to make the graphs either weight-balanced or
    unbalanced, as needed.}\label{fig:comm-topologies}
  %\vspace*{-3ex}
\end{figure}

\begin{figure}[htb]
  \centering
 \includegraphics[width=0.9\linewidth]{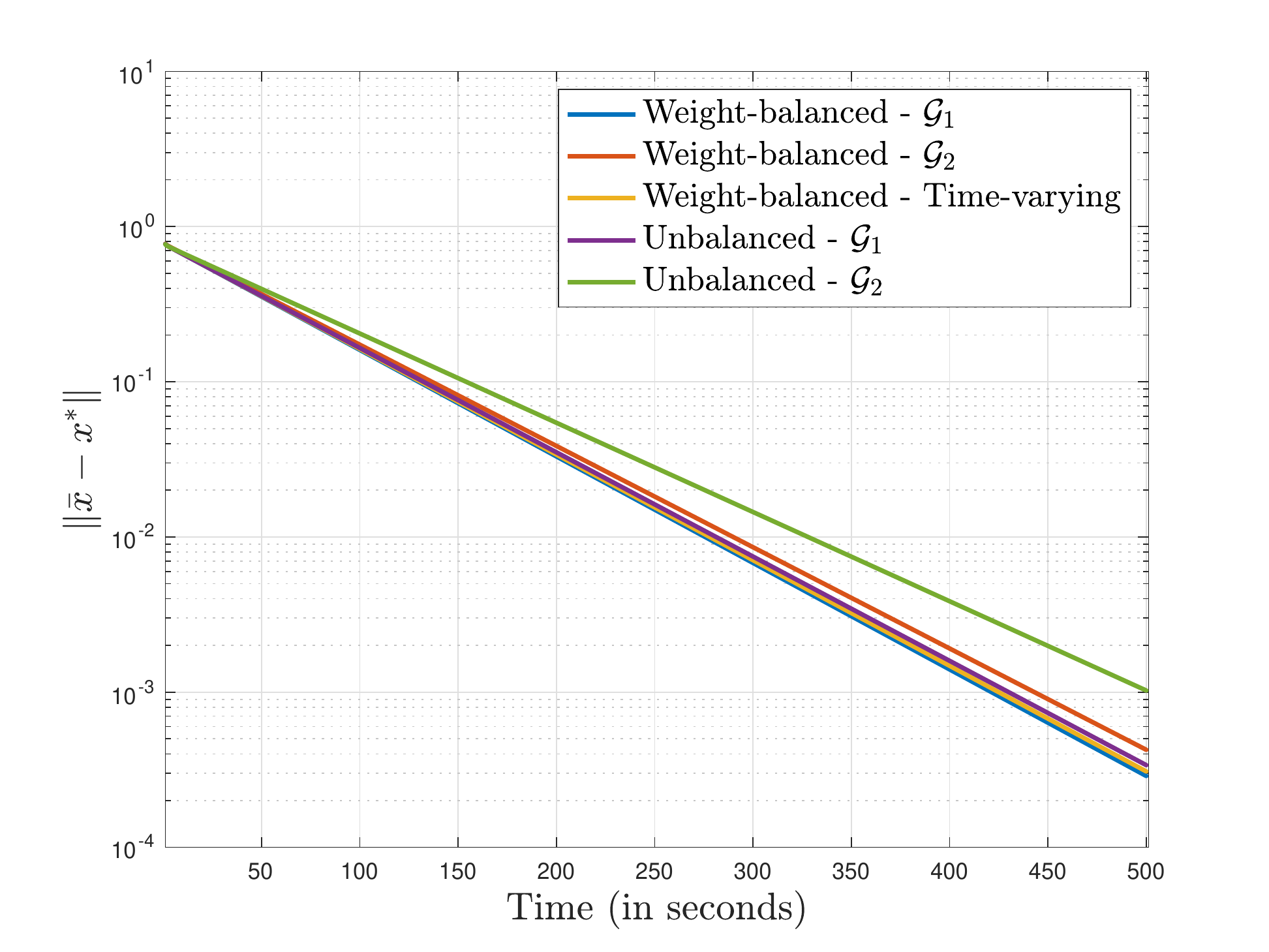}
 \caption{Evolution of the error between the actual solution and the
   average state using the proposed algorithms from initial condition
   $\xx(0)=\zero$, $\yy(0)=-\bb$, over the graphs shown in
   Fig.~\ref{fig:comm-topologies}.  The algorithms are implemented in
   discrete time with a stepsize of $2.5 \times 10^{-3}$, and the
   values of $\alpha=2$, $\beta=0.1$ and $\gamma=20$.  Straight lines
   correspond to exponential convergence.}
\label{fig:error}  
%\vspace*{-2ex}
\end{figure}

\section{Simulations}\label{sec:sims}
We consider 10 agents communicating over the digraphs shown in
Fig.~\ref{fig:comm-topologies}, seeking to solve
problem~\eqref{eq:problem} with $\{A_i\}_{i=1}^{10} \in \real^{5
  \times 5}$ and $\{b_i\}_{i=1}^{10} \in \real^5$.  Since the proposed
dynamics are in continuous time, we use a first-order Euler
discretization with stepsize $2.5 \times 10^{-3}$ for the MATLAB
implementation. The edge weights for various cases are adjusted to
make the graphs weight-balanced and unbalanced, resp.  For the
time-varying case, at every iteration, the communication graph is
switched randomly between $\G_1$ and $\G_2$.   
In Fig.~\ref{fig:error}, we plot the evolution of the error between
the actual solution of~\eqref{eq:problem} and the average state
$\bar{x} = \frac{1}{n} \onem^\top \xx$
using~\eqref{eq:distributed},~\eqref{eq:time_varying}
and~\eqref{eq:distributed_v}. The initial conditions for all the
algorithms are chosen according to Remark~\ref{remark:init}.  Even
though $\G_2$ (with $4.6$ as the minimum of the real parts of non-zero
eigenvalues of $\Lap$ and 
$\lambda_2(\Lap + \Lap^\top) = 7.6$, for the
weight-balanced case) is more connected than $\G_1$ (with $1.9$ as the
minimum of the real parts of non-zero eigenvalues of $\Lap$ and
$\lambda_2(\Lap + \Lap^\top) = 3.8$, for the weight-balanced case),
convergence is slower.  The error in the time-varying case is lower
and upper bounded by the error for $\G_1$ and $\G_2$, resp.

\section{Conclusions and Future Work}\label{sec:conclusions}
We have presented continuous-time algorithms to solve linear algebraic
equations whose problem data is represented as the summation of the
data of individual agents.  The proposed algorithms are distributed
over general directed networks, do not require the individual agent
matrices to be positive definite, and are guaranteed to converge to a
solution of the linear equation exponentially fast.  Future work will
involve formally characterizing the convergence when the null space 
condition is not satisfied, and
explore the design of distributed algorithms for finding least-square
solutions when exact ones do not exist, extension to cases where the
problem data is time-varying, and the communication graph is
unbalanced and time-varying.

\end{document}